\documentclass[12pt,a4paper,equation]{article}
\usepackage[greek,american]{babel}
\usepackage[iso-8859-7]{inputenc}
\usepackage{amssymb,amsfonts}
\usepackage[dvips]{graphicx}
\usepackage{amsmath}
\usepackage{epsfig}
\usepackage{latexsym}
\usepackage{makeidx}
\usepackage{pst-math,pst-xkey}

\numberwithin{equation}{section}

\makeindex
\begin{document}
\date{}
\author{{\bf Aristides V. Doumas}\\
%EndAName
Department of Mathematics\\
National Technical University of Athens\\
Zografou Campus\\
157 80 Athens, GREECE\\
\underline{aris.doumas@hotmail.com}\\
and\\
{\bf Vassilis G. Papanicolaou}\\
Department of Mathematics\\
National Technical University of Athens\\
Zografou Campus\\
157 80 Athens, GREECE\\
\underline{papanico@math.ntua.gr}}
\title{A Randomized Version of the Collatz $3x + 1$ Problem}
\maketitle
\begin{abstract}
We propose a stochastic version of the Collatz $3x + 1$ Problem.
\end{abstract}

\textbf{Keywords.} Collatz $3x + 1$ Problem; Markov chain.

\smallskip

\noindent\textbf{2010 AMS Mathematics Classification.} 60J10.

\section{The Collatz $3x + 1$ problem}

The classical Collatz $3x + 1$ Problem can be formulated as follows. Let $x$ be a positive odd integer. Consider the sequence
\begin{equation}
x_0 = x,
\qquad \qquad
x_n = \frac{3 x_{n-1} + 1}{2^{d_n}},
\quad
n \geq 1,
\label{A1}
\end{equation}
where $2^{d_n}$ is the highest power of $2$ dividing $3 x_{n-1} + 1$. Hence $\{x_n\}_{n=0}^{\infty}$ is a sequence of positive odd integers
(if, e.g., $x = 1$, then $x_n = 1$ for all $n$). Notice that $d_n \geq 1$ for all $n \geq 1$.

Suppose $L := \liminf x_n < \infty$. Then, since $x_n$ takes only positive integral values, we must have that $x_n = L$ for infinitely many
values of $n$. In particular, $x_k = x_{k+b} = L$ for some $k \geq 0$, $b \geq 1$. But, then, it follows from \eqref{A1} that
$x_{k+n} = x_{k+b+n}$ for all integers $n \geq 0$. Therefore, either
\begin{equation}
\lim_n x_n = \infty,
\label{A2}
\end{equation}
or the sequence $\{x_n\}_{n=0}^{\infty}$ is eventually periodic, namely there is a $b \geq 1$ and an $n_0 \geq 0$ such that
\begin{equation}
x_{n+b} = x_n
\qquad \text{for all }\; n \geq n_0.
\label{A3}
\end{equation}
Notice that, if $b=1$, i.e. if there is a $n_0$ such that $x_{n+1} = x_n$ for all $n \geq n_0$, then \eqref{A1} implies that
$(2^{d_{n+1}} - 3) x_n = 1$, which forces $x_n = 1$ for all $n \geq n_0$.

The Collatz $3x + 1$ Problem pertains to the behavior of the sequence $\{x_n\}_{n=0}^{\infty}$ as $n \to \infty$. One famous and longstanding open
question is whether there exists some initial value $x$ for which $\lim_n x_n = \infty$, while another open question is whether it is possible
to have an eventually periodic behavior with a (minimal) period $b > 1$.

The ultimate Collatz Conjecture is that, no matter what the initial value $x$ is, we always have that $x_n = 1$ for all $n$ sufficiently large.
Needless to say that the conjecture has beed verified for a huge set of initial values $x$.

\section{A randomized version of the problem}

Let $x$ be a positive odd integer. We consider the sequence
\begin{equation}
X_0 = x,
\qquad \qquad
X_n = \frac{3 X_{n-1} + \xi_n}{2^{d_n}},
\quad
n \geq 1,
\label{B1}
\end{equation}
where $\{\xi_n\}_{n=1}^{\infty}$ is a sequence of independent and identically distributed (i.i.d.) random variables taking odd integral values
$\geq -1$, and $2^{d_n}$ is the highest power of $2$ dividing $3 X_{n-1} + \xi_n$ (notice that we again have $d_n \geq 1$ for all $n \geq 1$).
Thus, $\{X_n\}_{n=0}^{\infty}$ is now a random sequence of positive odd integers.

Let us introduce the filtration
\begin{equation}
\mathcal{F}_0 := \{\emptyset, \Omega\},
\qquad \qquad
\mathcal{F}_n := \sigma(\xi_1, \dots, \xi_n),
\quad
n \geq 1,
\label{B3}
\end{equation}
where $(\Omega, \mathcal{F}, P)$ is the underlying probability space.
Clearly, by \eqref{B1} we have that the random variables $X_n$ and $d_n$ are $\mathcal{F}_n$-measurable for all $n \geq 1$. Notice that
\begin{equation}
\mathcal{F}^X_n := \sigma(X_1, \dots, X_n) \subset \mathcal{F}_n
\quad \text{and} \quad
\mathcal{F}^d_n := \sigma(d_1, \dots, d_n) \subset \mathcal{F}_n,
\quad
n \geq 1.
\label{B2a}
\end{equation}
Of course,
\begin{equation}
\mathcal{F}^X_n \vee \mathcal{F}^d_n = \mathcal{F}_n,
\qquad
n \geq 1,
\label{B2b}
\end{equation}
where $\mathcal{F}^X_n \vee \mathcal{F}^d_n$ denotes the $\sigma$-algebra generated by $\mathcal{F}^X_n$ and $\mathcal{F}^d_n$.

Formula \eqref{B1} implies that $\{X_n\}_{n=0}^{\infty}$ is a Markov chain with respect to $\mathcal{F}_n$,
whose state space is the set $\mathbb{N}_{\text{odd}}$ of positive odd integers. Actually, the two-dimensional process $\{(X_n, d_n)\}_{n=0}^{\infty}$ can be also viewed as a Markov chain with respect to $\mathcal{F}_n$ (the value of $d_0$ is irrelevant; furthermore, conditioning on $d_n$ is
irrelevant for $(X_{n+1}, d_{n+1})$).

The most natural case to examine first seems to be the choice $P\{\xi_n = -1\} = P\{\xi_n = 1\} = 1/2$
(or $P\{\xi_n = 1\} = P\{\xi_n = 3\} = 1/2$). Here, however, we will consider the rather easier case
\begin{equation}
P\{\xi_n = 1\} = P\{\xi_n = 3\} = P\{\xi_n = 5\} = P\{\xi_n = 7\} = \frac{1}{4}.
\label{B2}
\end{equation}

To begin our analysis, let us observe that for any positive odd integral value of $X_{n-1}$ we have
\begin{equation}
\{3 X_{n-1} + 1, \,\, 3 X_{n-1} + 3, \,\, 3 X_{n-1} + 5, \,\, 3 X_{n-1} + 7\} \equiv \{0, 2, 4, 6\} \ \; \text{mod}\; 8.
\label{B3a}
\end{equation}
Therefore, due to \eqref{B2} and the independence of $X_{n-1}$ and $\xi_n$ we have
\begin{equation}
P\{3 X_{n-1} + \xi_n \equiv k \ \text{mod}\; 8\} = \frac{1}{4}
\qquad \text{for }\;
k = 0, 2, 4, 6.
\label{B3b}
\end{equation}
The above formula motivates us to set
\begin{equation}
m_n := 3 \wedge d_n,
\qquad n \geq 1
\label{B4}
\end{equation}
(as usual, $a \wedge b$ denotes the minimum of $a$ and $b$), where $d_n$ is the random exponent appearing in \eqref{B1}. Then, \eqref{B1},
\eqref{B3b}, the Markov property of $(X_n, d_n)$, and the independence of $X_{n-1}$ and $\xi_n$ imply
\begin{equation}
P\{m_n = 1 \, | \, \mathcal{F}_{n-1}\} = P\{m_n = 1 \, | \, X_{n-1}\} =
P\{3 X_{n-1} + \, \xi_n \equiv 2\ \text{mod}\; 4 \, | \, X_{n-1}\} = \frac{1}{2}
\label{B5}
\end{equation}
for all $n \geq 1$. Likewise,
\begin{equation}
P\{m_n = 2 \, | \, \mathcal{F}_{n-1}\} = P\{m_n = 3 \, | \, \mathcal{F}_{n-1}\} = \frac{1}{4},
\qquad
n \geq 1.
\label{B6}
\end{equation}
Formulas \eqref{B5} and \eqref{B6} tell us that $m_n$ and $\mathcal{F}_{n-1}$ are independent for every $n \geq 1$. In particular
(since $m_n$ is $\mathcal{F}_n$-measurable for all $n \geq 1$) we have that $\{m_n\}_{n=1}^{\infty}$ is a sequence of
i.i.d. random variables with
\begin{equation}
P\{m_n = 1\} = \frac{1}{2},
\qquad
P\{m_n = 2\} = P\{m_n = 3\} = \frac{1}{4}.
\label{B7}
\end{equation}

\textbf{Proposition 1.} Let $\{X_n\}_{n=0}^{\infty}$ be the odd-integer-valued Markov chain introduced in \eqref{B1}, where
$\{\xi_n\}_{n=1}^{\infty}$ is a sequence of i.i.d. random variables whose common distribution is given by \eqref{B2}. Then,
\begin{equation}
\limsup X_n = \infty \ \; \text{a}.\text{s}.
\label{B8}
\end{equation}
for any initial value $X_0 = x$ (as usual, ``a.s." stands for ``almost surely", i.e. ``with probability $1$").

\smallskip

\textit{Proof}. Fix a constant $K > 0$ and then pick an integer $k \geq 1$ so that
\begin{equation}
\left(\frac{3}{2}\right)^k > K.
\label{B9}
\end{equation}
Since $\{m_n\}_{n=1}^{\infty}$ is a sequence of i.i.d. random variables whose common distribution is given by \eqref{B7},
an immediate consequence of the 2nd Borel-Cantelli Lemma is that
\begin{equation}
P\{m_n = m_{n+1} = \cdots = m_{n+k-1} = 1 \ \; \text{i.o.}\} = 1
\label{B10}
\end{equation}
(``i.o." stands for ``infinitely often", i.e. for infinitely many values of $n$).
But, if $m_n = m_{n+1} = \cdots = m_{n+k-1} = 1$, then, by \eqref{B1}, \eqref{B4}, \eqref{B9}, and the fact that ($X_0 = x \geq 1$ and)
$X_n, \xi_n \geq 1$ for every $n \geq 1$, we must also have that
\begin{equation}
X_{n+k} > \left(\frac{3}{2}\right)^k X_n > K.
\label{B11}
\end{equation}
Therefore, \eqref{B10} implies that
\begin{equation}
P\{X_{n+k} > K \ \; \text{i.o.}\} = 1.
\label{B12}
\end{equation}
From formula \eqref{B12} we get
\begin{equation*}
\limsup X_n \geq K \ \; \text{a.s.}
%\label{B13}
\end{equation*}
and since $K$ is arbitrary, the proposition follows from the fact that
\begin{equation*}
\{\limsup X_n = \infty\}  = \bigcap_{K=1}^{\infty} \{\limsup X_n \geq K\}.
%\label{B14}
\end{equation*}
\hfill $\blacksquare$

\smallskip

\textbf{Remark 1.} In the case $P\{\xi_n = -1\} = P\{\xi_n = 1\} = 1/2$ things are quite different since now $1$ is an absorbing (or trapping) state,
i.e. if $X_n = 1$ for some $n$, then $X_{n+k} = 1$ for all $k \geq 0$. Hence, \eqref{B11} and, consequently, \eqref{B8} are not valid. In fact, since for any initial state $x$ there is always a positive probability that $X_n = 1$ for some $n$, it follows that $P\{\limsup X_n = \infty\} < 1$. A
natural \textbf{open question} here, in the spirit of the Collatz Problem, is whether $P\{\limsup X_n = \infty\} = 0$.

\smallskip

Let us now continue the analysis of the case \eqref{B2}. Formula \eqref{B1} can be written as
\begin{equation}
X_n = \frac{3}{2^{d_n}} \left(1 + \frac{\xi_n}{3 X_{n-1}}\right)  X_{n-1}.
\label{C9}
\end{equation}
Thus, in view of \eqref{B2} and \eqref{B4} we have
\begin{equation}
1 \leq X_n \leq \frac{3}{2^{m_n}} \left(1 + \frac{7}{3 X_{n-1}}\right) X_{n-1}.
\label{C10}
\end{equation}
Due to the multiplicative form of the formulas it is convenient to set
\begin{equation}
Y_n := \ln X_n,
\qquad
n \geq 0.
\label{C11}
\end{equation}
Then, inequality \eqref{C10} is equivalent to
\begin{equation}
0 \leq Y_n \leq Y_{n-1} + \ln\left(1 + \frac{7}{3} e^{-Y_{n-1}}\right) + \ln3 - m_n \ln2.
\label{C12}
\end{equation}
By \eqref{B7} we have $E[m_n] = 7/4$ (and $V[m_n] = 11/16$). Thus
\begin{equation}
E[\ln3 - m_n \ln2] = \ln3 - \frac{7 \ln2}{4} = -\frac{1}{4} \ln\left(\frac{128}{81}\right) \simeq -0.1144.
\label{C13}
\end{equation}
Now, let us fix an $\varepsilon$ such that $0 < \varepsilon < -E[\ln3 - m_n \ln2]$. Then, it is easy to see that there is an $M > 1$ such that
$\ln(1 + (7/3) e^{-y}) \leq \varepsilon$ for all $y \geq M$. For example, if $\varepsilon = 1/10$, then it suffices to take $M = \ln23$.
With such values of $\varepsilon$ and $M$, formula \eqref{C12} implies
\begin{equation}
0 \leq Y_n \leq Y_{n-1} + \varepsilon + \ln3 - m_n \ln2,
\qquad \text{if }\; Y_{n-1} \geq M,
\label{C14a}
\end{equation}
and
\begin{equation}
0 \leq Y_n \leq Y_{n-1} + \ln5,
\qquad \text{if }\; Y_{n-1} < M.
\label{C14b}
\end{equation}
For notational convenience we prefer to write formula \eqref{C14a} in the form
\begin{equation}
0 \leq Y_n \leq Y_{n-1} + W_n,
\qquad \text{if }\; Y_{n-1} \geq M,
\label{C14aa}
\end{equation}
where
\begin{equation}
W_n := \varepsilon + \ln3 - m_n \ln2,
\qquad
n \geq 0.
\label{C14aaa}
\end{equation}
From the properties of $\{m_n\}_{n=1}^{\infty}$ it follows immediately that $\{W_n\}_{n=1}^{\infty}$ is a sequence of i.i.d. random variables with
\begin{equation}
\mu := E[W_n] = \varepsilon + E[\ln3 - m_n \ln2] < 0
\label{C15}
\end{equation}
and, furthermore, that $Y_{n-1}$ and $W_n$ are independent for every $n \geq 1$.

The following proposition is in the spirit of the Collatz Conjecture.

\smallskip

\textbf{Proposition 2.} Let $\{X_n\}_{n=0}^{\infty}$ be the Markov chain of Proposition 1. Then, the state $1$ is positive recurrent \cite{D}.
In particular
\begin{equation}
P\{X_n = 1 \ \; \text{i.o.}\} = 1
\label{C16}
\end{equation}
for any initial value $X_0 = x$.

\smallskip

\textit{Proof}. For an $\varepsilon$ and an $M$ as above let
\begin{equation}
N_1 := \inf\{n \geq 0 \, : \, Y_n \geq M\}
\qquad \text{and} \qquad
D_1 := \inf\{n > N_1 \, : \, Y_n < M\}
\label{C17}
\end{equation}
(if $Y_0 = \ln x \geq M$, then $N_1 = 0$). Notice that $N_1$ and $D_1$ are stopping times of the Markov chain $\{Y_n\}_{n=0}^{\infty}$ and,
hence, of the filtration $\{\mathcal{F}_n\}_{n=1}^{\infty}$. By Proposition 1 we have that
\begin{equation}
N_1 < \infty \ \; \text{a}.\text{s}.
\label{C18}
\end{equation}
Actually, much more is true. Since for any fixed $k \geq 1$ we have that
\begin{equation*}
P\{m_n = m_{n+1} = \cdots = m_{n+k-1} = 1\} > 0,
%\label{C18a}
\end{equation*}
it follows that there is an $\alpha_0 > 0$ such that
\begin{equation}
E[e^{\alpha N_1}] < \infty
\qquad \text{for all }\;
\alpha < \alpha_0.
\label{C18b}
\end{equation}
In particular,
\begin{equation}
E[N_1] < \infty.
\label{C18c}
\end{equation}
Suppose now that $\omega \in \{D_1 = \infty\}$. Then, for such $\omega$'s we must have $Y_{N_1 + n} \geq M$ for all $n \geq 0$ and, hence,
formula \eqref{C14aa} becomes
\begin{equation*}
Y_{N_1 + n+1} - Y_{N_1 + n} \leq  W_{N_1 + n+1},
\qquad \text{for all }\; n \geq 0
%\label{C19a}
\end{equation*}
or
\begin{equation}
Y_{N_1 + n} - Y_{N_1} \leq  \sum_{j=1}^n W_{N_1 + j},
\qquad \text{for all }\; n \geq 1.
\label{C19}
\end{equation}
However, $\{W_{N_1 + j}\}_{j=1}^{\infty}$ is a sequence of i.i.d. random variables whose common distribution is that of $W_n$ \cite{D}.
In particular, the common expectation of $W_{N_1 + j}$, $j \geq 1$, is strictly negative. As a consequence of these facts we have that the
event of \eqref{C19} has probability $0$ and, furthermore, there is a $\beta_0 > 0$ such that
\begin{equation}
E[e^{\beta D_1}] < \infty
\qquad \text{for all }\;
\beta < \beta_0.
\label{C20}
\end{equation}
In particular,
\begin{equation}
E[D_1] < \infty
\label{C20a}
\end{equation}
and
\begin{equation}
D_1 < \infty \ \; \text{a}.\text{s}.
\label{C20b}
\end{equation}
We can then introduce the stopping times
\begin{equation}
N_k := \inf\{n > D_{k-1} \, : \, Y_n \geq M\},
\qquad
D_k := \inf\{n > N_k \, : \, Y_n < M\},
\quad
k \geq 2.
\label{C21}
\end{equation}
As in the case of $N_1$ and $D_1$, we again have that, there are $\alpha_0 > 0$ and $\beta_0 > 0$ independent of $k$ such that,
for all $k \geq 2$ we have
\begin{equation}
E[e^{\alpha N_k}] < \infty,
\quad
E[e^{\beta D_k}] < \infty
\qquad \text{for every }\;
\alpha < \alpha_0, \ \; \beta < \beta_0.
\label{C22}
\end{equation}
In particular,
\begin{equation}
E[N_k] < \infty,
\qquad \qquad
E[D_k] < \infty
\label{C20a}
\end{equation}
and
\begin{equation}
N_k < \infty, \ D_k < \infty \ \; \text{a}.\text{s}.
\qquad
\label{C22b}
\end{equation}

Let $\mathcal{O}_M$ denote the set of odd positive integers which are less than $e^M$ (definitely $1 \in \mathcal{O}_M$). Then, the
above analysis implies
\begin{equation}
P\{X_n \in \mathcal{O}_M \ \; \text{i.o.}\} = 1.
\label{C23}
\end{equation}
Since $\mathcal{O}_M$ is a finite set, some state $r$ in $\mathcal{O}_M$ must be recurrent; actually positive recurrent due to \eqref{C20a}
\cite{D}. But, then, since there is a nonzero probability for the Markov chain $X_n$ to go from any $r \in \mathcal{O}_M$ to $1$,
it follows that $1$ is positive recurrent \cite{D}.
\hfill $\blacksquare $

\smallskip

The next proposition generalizes Proposition 2.

\smallskip

\textbf{Proposition 3.} Let $\{X_n\}_{n=0}^{\infty}$ be the Markov chain of Propositions 1 and 2. Then,
all states in $\mathbb{N}_{\text{odd}}$ are positive recurrent.

\smallskip

\textit{Proof}. By Proposition 2 we can assume without loss of generality that $X_0 = 1$. Then, we need to show that all states in
$\mathbb{N}_{\text{odd}}$ can be reached by $\{X_n\}_{n=0}^{\infty}$ with nonzero probability.

Let $m$ be the smallest odd integer which cannot be reached, namely $P\{X_n = m\} = 0$ for all $n \geq 1$. By Proposition 2 we have
that $m \geq 3$. Therefore, we should have one of the following three possibilities:
\begin{equation*}
m = 6k + 3 \quad \text{or} \quad m = 6k + 5 \quad \text{or} \quad m = 6k + 7 \quad \qquad \text{for some }\; k \geq 0.
%\label{A3}
\end{equation*}
(i) Suppose $m = 6k +3$ for some $k \geq 0$. Then, $4k + 1 < m$ and since $4k + 1$ is odd, we must have that $P\{X_n = 4k + 1\} > 0$
for some $n$. But, then, since
\begin{equation*}
3(4k + 1) + 3 = 12k + 6 = 2(6k + 3) = 2m,
%\label{A4}
\end{equation*}
we should have by \eqref{B1} that $P\{X_{n+1} = m\} > 0$, a contradiction.

(ii) Next, suppose $m = 6k +5$ for some $k \geq 0$. Then, $4k + 3 < m$ and since $4k + 3$ is odd, we must have that $P\{X_n = 4k + 3\} > 0$
for some $n$. But, then, since
\begin{equation*}
3(4k + 3) + 1 = 12k + 10 = 2(6k + 5) = 2m,
%\label{A5}
\end{equation*}
we should have by \eqref{B1} that $P\{X_{n+1} = m\} > 0$, a contradiction.

(iii) Finally, suppose $m = 6k +7$ for some $k \geq 0$. Then, $4k + 3 < m$ and since $4k + 3$ is odd, we must have that  $P\{X_n = 4k + 3\} > 0$
for some $n$. But, then, since
\begin{equation*}
3(4k + 3) + 5 = 12k + 14 = 2(6k + 7) = 2m,
%\label{A6}
\end{equation*}
we should have by \eqref{B1} that $P\{X_{n+1} = m\} > 0$, again a contradiction.

Therefore, in all three possibilities for $m$ we have reached a contradiction. It follows that such an $m$ cannot exist.
\hfill $\blacksquare $

\textbf{Remark 2.} A side result of Proposition 3 is that the Markov chain $\{X_n\}_{n=0}^{\infty}$ is irreducible \cite{D}. In the case
where $\xi_n$ takes the values $1$, $3$, and $5$ only with positive probabilities, given $X_0 = 1$, one can use the idea of the proof of
Proposition 3 in order to show that the associated Markov chain is irreducible. However, it is an \textbf{open question} whether Proposition 2
is valid in that case. As for the case where $\xi_n$ takes only the values $1$ and $3$ with positive probabilities, given $X_0 = 1$, even the
irreducibility is an \textbf{open question}.

\smallskip

\textbf{Final Comments.} Since $P\{X_{n+1} = 1 \, | \, X_n = 1\} = 1/2 > 0$, it follows from Proposition 3 that the Markov chain
$\{X_n\}_{n=0}^{\infty}$ is aperiodic \cite{D}. Also, again by Proposition 3 we have that $\{X_n\}_{n=0}^{\infty}$ has a stationary distribution
$\pi$ \cite{D}. Finally, the existence of $\pi$ together with the aperiodicity (and the irreducibility mentioned in Remark 2) imply \cite{D} that
\begin{equation*}
P\{X_n = y \, | \, X_0 = x\} \to \pi(y)
\qquad \text{as }\; n \to \infty,
%\label{A7}
\end{equation*}
for every $x, y \in \mathbb{N}_{\text{odd}}$.

\end{document}